
\input epsf




\overfullrule=0pt
\magnification=1200
\hsize=11.25cm    
\vsize=18cm
\hoffset=1cm
\font\grand=cmr10 at 14pt
\font\petit=cmr10 at 8pt

\def\N{\noindent}
\def\S{\smallskip \par}
\def\M{\medskip \par}
\def\B{\bigskip \par}
\def\BB{\bigskip\bigskip\par}

\def\oo{\omega}
\def\DD{\Delta}
\def\mm{{}^{-1}}
\def\ss{\sigma}

\def\ee{\epsilon}


\def\RR{{\bf R}}

\def\RR{{\bf R}}

\def\KKK{{\cal K}}

\def\sqr#1#2{{\vcenter{\vbox{\hrule height.#2pt
\hbox{\vrule width .#2pt height#1pt \kern#1pt
\vrule width.#2pt}
\hrule height.#2pt}}}}
\def\square{\mathchoice\sqr64\sqr64\sqr{4.2}3\sqr33}

\def\proj{\mathop{\rightarrow \!\!\!\!\!\!\! \rightarrow}}


\def\row#1#2#3{(#1_{#2},\ldots,#1_{#3})}

\catcode`@=11
\def\@height{height}
\def\@depth{depth}
\def\@width{width}

\newcount\@tempcnta
\newcount\@tempcntb

\newdimen\@tempdima
\newdimen\@tempdimb

\newbox\@tempboxa

\def\@ifnextchar#1#2#3{\let\@tempe #1\def\@tempa{#2}\def\@tempb{#3}\futurelet
    \@tempc\@ifnch}
\def\@ifnch{\ifx \@tempc \@sptoken \let\@tempd\@xifnch
      \else \ifx \@tempc \@tempe\let\@tempd\@tempa\else\let\@tempd\@tempb\fi
      \fi \@tempd}
\def\@ifstar#1#2{\@ifnextchar *{\def\@tempa*{#1}\@tempa}{#2}}

\def\@whilenoop#1{}
\def\@whilenum#1\do #2{\ifnum #1\relax #2\relax\@iwhilenum{#1\relax 
     #2\relax}\fi}
\def\@iwhilenum#1{\ifnum #1\let\@nextwhile=\@iwhilenum 
         \else\let\@nextwhile=\@whilenoop\fi\@nextwhile{#1}}

\def\@whiledim#1\do #2{\ifdim #1\relax#2\@iwhiledim{#1\relax#2}\fi}
\def\@iwhiledim#1{\ifdim #1\let\@nextwhile=\@iwhiledim 
        \else\let\@nextwhile=\@whilenoop\fi\@nextwhile{#1}}

\newdimen\@wholewidth
\newdimen\@halfwidth
\newdimen\unitlength \unitlength =1pt
\newbox\@picbox
\newdimen\@picht

\def\@nnil{\@nil}
\def\@empty{}
\def\@fornoop#1\@@#2#3{}

\def\@for#1:=#2\do#3{\edef\@fortmp{#2}\ifx\@fortmp\@empty \else
    \expandafter\@forloop#2,\@nil,\@nil\@@#1{#3}\fi}

\def\@forloop#1,#2,#3\@@#4#5{\def#4{#1}\ifx #4\@nnil \else
       #5\def#4{#2}\ifx #4\@nnil \else#5\@iforloop #3\@@#4{#5}\fi\fi}

\def\@iforloop#1,#2\@@#3#4{\def#3{#1}\ifx #3\@nnil 
       \let\@nextwhile=\@fornoop \else
      #4\relax\let\@nextwhile=\@iforloop\fi\@nextwhile#2\@@#3{#4}}

\def\@tfor#1:=#2\do#3{\xdef\@fortmp{#2}\ifx\@fortmp\@empty \else
    \@tforloop#2\@nil\@nil\@@#1{#3}\fi}
\def\@tforloop#1#2\@@#3#4{\def#3{#1}\ifx #3\@nnil 
       \let\@nextwhile=\@fornoop \else
      #4\relax\let\@nextwhile=\@tforloop\fi\@nextwhile#2\@@#3{#4}}

\def\@makepicbox(#1,#2){\leavevmode\@ifnextchar 
   [{\@imakepicbox(#1,#2)}{\@imakepicbox(#1,#2)[]}}

\long\def\@imakepicbox(#1,#2)[#3]#4{\vbox to#2\unitlength
   {\let\mb@b\vss \let\mb@l\hss\let\mb@r\hss
    \let\mb@t\vss
    \@tfor\@tempa :=#3\do{\expandafter\let
        \csname mb@\@tempa\endcsname\relax}%
\mb@t\hbox to #1\unitlength{\mb@l #4\mb@r}\mb@b}}

\def\picture(#1,#2){\@ifnextchar({\@picture(#1,#2)}{\@picture(#1,#2)(0,0)}}

\def\@picture(#1,#2)(#3,#4){\@picht #2\unitlength
\setbox\@picbox\hbox to #1\unitlength\bgroup 
\hskip -#3\unitlength \lower #4\unitlength \hbox\bgroup\ignorespaces}

\def\endpicture{\egroup\hss\egroup\ht\@picbox\@picht
\dp\@picbox\z@\leavevmode\box\@picbox}

\long\def\put(#1,#2)#3{\@killglue\raise#2\unitlength\hbox to \z@{\kern
#1\unitlength #3\hss}\ignorespaces}

\long\def\multiput(#1,#2)(#3,#4)#5#6{\@killglue\@multicnt=#5\relax
\@xdim=#1\unitlength
\@ydim=#2\unitlength
\@whilenum \@multicnt > 0\do
{\raise\@ydim\hbox to \z@{\kern
\@xdim #6\hss}\advance\@multicnt \m@ne\advance\@xdim
#3\unitlength\advance\@ydim #4\unitlength}\ignorespaces}

\def\@killglue{\unskip\@whiledim \lastskip >\z@\do{\unskip}}

\def\thinlines{\let\@linefnt\tenln \let\@circlefnt\tencirc
  \@wholewidth\fontdimen8\tenln \@halfwidth .5\@wholewidth}
\def\thicklines{\let\@linefnt\tenlnw \let\@circlefnt\tencircw
  \@wholewidth\fontdimen8\tenlnw \@halfwidth .5\@wholewidth}

\def\linethickness#1{\@wholewidth #1\relax \@halfwidth .5\@wholewidth}

\def\shortstack{\@ifnextchar[{\@shortstack}{\@shortstack[c]}}

\def\@shortstack[#1]{\leavevmode
\vbox\bgroup\baselineskip-1pt\lineskip 3pt\let\mb@l\hss
\let\mb@r\hss \expandafter\let\csname mb@#1\endcsname\relax
\let\\\@stackcr\@ishortstack}

\def\@ishortstack#1{\halign{\mb@l ##\unskip\mb@r\cr #1\crcr}\egroup}

\def\@stackcr{\@ifstar{\@ixstackcr}{\@ixstackcr}}
\def\@ixstackcr{\@ifnextchar[{\@istackcr}{\cr\ignorespaces}}

\def\@istackcr[#1]{\cr\noalign{\vskip #1}\ignorespaces}

\newif\if@negarg

\def\droite(#1,#2)#3{\@xarg #1\relax \@yarg #2\relax
\@linelen=#3\unitlength
\ifnum\@xarg =0 \@vline 
  \else \ifnum\@yarg =0 \@hline \else \@sline\fi
\fi}

\def\@sline{\ifnum\@xarg< 0 \@negargtrue \@xarg -\@xarg \@yyarg -\@yarg
  \else \@negargfalse \@yyarg \@yarg \fi
\ifnum \@yyarg >0 \@tempcnta\@yyarg \else \@tempcnta -\@yyarg \fi
\ifnum\@tempcnta>6 \@badlinearg\@tempcnta0 \fi
\ifnum\@xarg>6 \@badlinearg\@xarg 1 \fi
\setbox\@linechar\hbox{\@linefnt\@getlinechar(\@xarg,\@yyarg)}%
\ifnum \@yarg >0 \let\@upordown\raise \@clnht\z@
   \else\let\@upordown\lower \@clnht \ht\@linechar\fi
\@clnwd=\wd\@linechar
\if@negarg \hskip -\wd\@linechar \def\@tempa{\hskip -2\wd\@linechar}\else
     \let\@tempa\relax \fi
\@whiledim \@clnwd <\@linelen \do
  {\@upordown\@clnht\copy\@linechar
   \@tempa
   \advance\@clnht \ht\@linechar
   \advance\@clnwd \wd\@linechar}%
\advance\@clnht -\ht\@linechar
\advance\@clnwd -\wd\@linechar
\@tempdima\@linelen\advance\@tempdima -\@clnwd
\@tempdimb\@tempdima\advance\@tempdimb -\wd\@linechar
\if@negarg \hskip -\@tempdimb \else \hskip \@tempdimb \fi
\multiply\@tempdima \@m
\@tempcnta \@tempdima \@tempdima \wd\@linechar \divide\@tempcnta \@tempdima
\@tempdima \ht\@linechar \multiply\@tempdima \@tempcnta
\divide\@tempdima \@m
\advance\@clnht \@tempdima
\ifdim \@linelen <\wd\@linechar
   \hskip \wd\@linechar
  \else\@upordown\@clnht\copy\@linechar\fi}

\def\@hline{\ifnum \@xarg <0 \hskip -\@linelen \fi
\vrule \@height \@halfwidth \@depth \@halfwidth \@width \@linelen
\ifnum \@xarg <0 \hskip -\@linelen \fi}

\def\@getlinechar(#1,#2){\@tempcnta#1\relax\multiply\@tempcnta 8
\advance\@tempcnta -9 \ifnum #2>0 \advance\@tempcnta #2\relax\else
\advance\@tempcnta -#2\relax\advance\@tempcnta 64 \fi
\char\@tempcnta}

\def\vector(#1,#2)#3{\@xarg #1\relax \@yarg #2\relax
\@tempcnta \ifnum\@xarg<0 -\@xarg\else\@xarg\fi
\ifnum\@tempcnta<5\relax
\@linelen=#3\unitlength
\ifnum\@xarg =0 \@vvector 
  \else \ifnum\@yarg =0 \@hvector \else \@svector\fi
\fi
\else\@badlinearg\fi}

\def\@hvector{\@hline\hbox to 0pt{\@linefnt 
\ifnum \@xarg <0 \@getlarrow(1,0)\hss\else
    \hss\@getrarrow(1,0)\fi}}

\def\@vvector{\ifnum \@yarg <0 \@downvector \else \@upvector \fi}

\def\@svector{\@sline
\@tempcnta\@yarg \ifnum\@tempcnta <0 \@tempcnta=-\@tempcnta\fi
\ifnum\@tempcnta <5
  \hskip -\wd\@linechar
  \@upordown\@clnht \hbox{\@linefnt  \if@negarg 
  \@getlarrow(\@xarg,\@yyarg) \else \@getrarrow(\@xarg,\@yyarg) \fi}%
\else\@badlinearg\fi}

\def\@getlarrow(#1,#2){\ifnum #2 =\z@ \@tempcnta='33\else
\@tempcnta=#1\relax\multiply\@tempcnta \sixt@@n \advance\@tempcnta
-9 \@tempcntb=#2\relax\multiply\@tempcntb \tw@
\ifnum \@tempcntb >0 \advance\@tempcnta \@tempcntb\relax
\else\advance\@tempcnta -\@tempcntb\advance\@tempcnta 64
\fi\fi\char\@tempcnta}

\def\@getrarrow(#1,#2){\@tempcntb=#2\relax
\ifnum\@tempcntb < 0 \@tempcntb=-\@tempcntb\relax\fi
\ifcase \@tempcntb\relax \@tempcnta='55 \or 
\ifnum #1<3 \@tempcnta=#1\relax\multiply\@tempcnta
24 \advance\@tempcnta -6 \else \ifnum #1=3 \@tempcnta=49
\else\@tempcnta=58 \fi\fi\or 
\ifnum #1<3 \@tempcnta=#1\relax\multiply\@tempcnta
24 \advance\@tempcnta -3 \else \@tempcnta=51\fi\or 
\@tempcnta=#1\relax\multiply\@tempcnta
\sixt@@n \advance\@tempcnta -\tw@ \else
\@tempcnta=#1\relax\multiply\@tempcnta
\sixt@@n \advance\@tempcnta 7 \fi\ifnum #2<0 \advance\@tempcnta 64 \fi
\char\@tempcnta}

\def\@vline{\ifnum \@yarg <0 \@downline \else \@upline\fi}

\def\@upline{\hbox to \z@{\hskip -\@halfwidth \vrule \@width \@wholewidth
   \@height \@linelen \@depth \z@\hss}}

\def\@downline{\hbox to \z@{\hskip -\@halfwidth \vrule \@width \@wholewidth
   \@height \z@ \@depth \@linelen \hss}}

\def\@upvector{\@upline\setbox\@tempboxa\hbox{\@linefnt\char'66}\raise 
     \@linelen \hbox to\z@{\lower \ht\@tempboxa\box\@tempboxa\hss}}

\def\@downvector{\@downline\lower \@linelen
      \hbox to \z@{\@linefnt\char'77\hss}}

\def\dashbox#1(#2,#3){\leavevmode\hbox to \z@{\baselineskip \z@%
\lineskip \z@%
\@dashdim=#2\unitlength%
\@dashcnt=\@dashdim \advance\@dashcnt 200
\@dashdim=#1\unitlength\divide\@dashcnt \@dashdim
\ifodd\@dashcnt\@dashdim=\z@%
\advance\@dashcnt \@ne \divide\@dashcnt \tw@ 
\else \divide\@dashdim \tw@ \divide\@dashcnt \tw@
\advance\@dashcnt \m@ne
\setbox\@dashbox=\hbox{\vrule \@height \@halfwidth \@depth \@halfwidth
\@width \@dashdim}\put(0,0){\copy\@dashbox}%
\put(0,#3){\copy\@dashbox}%
\put(#2,0){\hskip-\@dashdim\copy\@dashbox}%
\put(#2,#3){\hskip-\@dashdim\box\@dashbox}%
\multiply\@dashdim 3 
\fi
\setbox\@dashbox=\hbox{\vrule \@height \@halfwidth \@depth \@halfwidth
\@width #1\unitlength\hskip #1\unitlength}\@tempcnta=0
\put(0,0){\hskip\@dashdim \@whilenum \@tempcnta <\@dashcnt
\do{\copy\@dashbox\advance\@tempcnta \@ne }}\@tempcnta=0
\put(0,#3){\hskip\@dashdim \@whilenum \@tempcnta <\@dashcnt
\do{\copy\@dashbox\advance\@tempcnta \@ne }}%
\@dashdim=#3\unitlength%
\@dashcnt=\@dashdim \advance\@dashcnt 200
\@dashdim=#1\unitlength\divide\@dashcnt \@dashdim
\ifodd\@dashcnt \@dashdim=\z@%
\advance\@dashcnt \@ne \divide\@dashcnt \tw@
\else
\divide\@dashdim \tw@ \divide\@dashcnt \tw@
\advance\@dashcnt \m@ne
\setbox\@dashbox\hbox{\hskip -\@halfwidth
\vrule \@width \@wholewidth 
\@height \@dashdim}\put(0,0){\copy\@dashbox}%
\put(#2,0){\copy\@dashbox}%
\put(0,#3){\lower\@dashdim\copy\@dashbox}%
\put(#2,#3){\lower\@dashdim\copy\@dashbox}%
\multiply\@dashdim 3
\fi
\setbox\@dashbox\hbox{\vrule \@width \@wholewidth 
\@height #1\unitlength}\@tempcnta0
\put(0,0){\hskip -\@halfwidth \vbox{\@whilenum \@tempcnta < \@dashcnt
\do{\vskip #1\unitlength\copy\@dashbox\advance\@tempcnta \@ne }%
\vskip\@dashdim}}\@tempcnta0
\put(#2,0){\hskip -\@halfwidth \vbox{\@whilenum \@tempcnta< \@dashcnt
\relax\do{\vskip #1\unitlength\copy\@dashbox\advance\@tempcnta \@ne }%
\vskip\@dashdim}}}\@makepicbox(#2,#3)}

\newif\if@ovt 
\newif\if@ovb 
\newif\if@ovl 
\newif\if@ovr 
\newdimen\@ovxx
\newdimen\@ovyy
\newdimen\@ovdx
\newdimen\@ovdy
\newdimen\@ovro
\newdimen\@ovri

\def\@getcirc#1{\@tempdima #1\relax \advance\@tempdima 2pt\relax
  \@tempcnta\@tempdima
  \@tempdima 4pt\relax \divide\@tempcnta\@tempdima
  \ifnum \@tempcnta > 10\relax \@tempcnta 10\relax\fi
  \ifnum \@tempcnta >\z@ \advance\@tempcnta\m@ne
    \else \@warning{Oval too small}\fi
  \multiply\@tempcnta 4\relax
  \setbox \@tempboxa \hbox{\@circlefnt
  \char \@tempcnta}\@tempdima \wd \@tempboxa}

\def\@put#1#2#3{\raise #2\hbox to \z@{\hskip #1#3\hss}}

\def\oval(#1,#2){\@ifnextchar[{\@oval(#1,#2)}{\@oval(#1,#2)[]}}

\def\@oval(#1,#2)[#3]{\begingroup\boxmaxdepth \maxdimen
  \@ovttrue \@ovbtrue \@ovltrue \@ovrtrue
  \@tfor\@tempa :=#3\do{\csname @ov\@tempa false\endcsname}\@ovxx
  #1\unitlength \@ovyy #2\unitlength
  \@tempdimb \ifdim \@ovyy >\@ovxx \@ovxx\else \@ovyy \fi
  \advance \@tempdimb -2pt\relax  
  \@getcirc \@tempdimb
  \@ovro \ht\@tempboxa \@ovri \dp\@tempboxa
  \@ovdx\@ovxx \advance\@ovdx -\@tempdima \divide\@ovdx \tw@
  \@ovdy\@ovyy \advance\@ovdy -\@tempdima \divide\@ovdy \tw@
  \@circlefnt \setbox\@tempboxa
  \hbox{\if@ovr \@ovvert32\kern -\@tempdima \fi
  \if@ovl \kern \@ovxx \@ovvert01\kern -\@tempdima \kern -\@ovxx \fi
  \if@ovt \@ovhorz \kern -\@ovxx \fi
  \if@ovb \raise \@ovyy \@ovhorz \fi}\advance\@ovdx\@ovro
  \advance\@ovdy\@ovro \ht\@tempboxa\z@ \dp\@tempboxa\z@
  \@put{-\@ovdx}{-\@ovdy}{\box\@tempboxa}%
  \endgroup}

\def\@ovvert#1#2{\vbox to \@ovyy{%
    \if@ovb \@tempcntb \@tempcnta \advance \@tempcntb by #1\relax
      \kern -\@ovro \hbox{\char \@tempcntb}\nointerlineskip
    \else \kern \@ovri \kern \@ovdy \fi
    \leaders\vrule width \@wholewidth\vfil \nointerlineskip
    \if@ovt \@tempcntb \@tempcnta \advance \@tempcntb by #2\relax
      \hbox{\char \@tempcntb}%
    \else \kern \@ovdy \kern \@ovro \fi}}

\def\@ovhorz{\hbox to \@ovxx{\kern \@ovro
    \if@ovr \else \kern \@ovdx \fi
    \leaders \hrule height \@wholewidth \hfil
    \if@ovl \else \kern \@ovdx \fi
    \kern \@ovri}}

\def\circle{\@ifstar{\@dot}{\@circle}}
\def\@circle#1{\begingroup \boxmaxdepth \maxdimen \@tempdimb #1\unitlength
   \ifdim \@tempdimb >15.5pt\relax \@getcirc\@tempdimb
      \@ovro\ht\@tempboxa 
     \setbox\@tempboxa\hbox{\@circlefnt
      \advance\@tempcnta\tw@ \char \@tempcnta
      \advance\@tempcnta\m@ne \char \@tempcnta \kern -2\@tempdima
      \advance\@tempcnta\tw@
      \raise \@tempdima \hbox{\char\@tempcnta}\raise \@tempdima
        \box\@tempboxa}\ht\@tempboxa\z@ \dp\@tempboxa\z@
      \@put{-\@ovro}{-\@ovro}{\box\@tempboxa}%
   \else  \@circ\@tempdimb{96}\fi\endgroup}

\def\@dot#1{\@tempdimb #1\unitlength \@circ\@tempdimb{112}}

\def\@circ#1#2{\@tempdima #1\relax \advance\@tempdima .5pt\relax
   \@tempcnta\@tempdima \@tempdima 1pt\relax
   \divide\@tempcnta\@tempdima 
   \ifnum\@tempcnta > 15\relax \@tempcnta 15\relax \fi    
   \ifnum \@tempcnta >\z@ \advance\@tempcnta\m@ne\fi
   \advance\@tempcnta #2\relax
   \@circlefnt \char\@tempcnta}

\font\tenln line10
\font\tencirc lcircle10
\font\tenlnw linew10
\font\tencircw lcirclew10

\thinlines   

\newcount\@xarg
\newcount\@yarg
\newcount\@yyarg
\newcount\@multicnt 
\newdimen\@xdim
\newdimen\@ydim
\newbox\@linechar
\newdimen\@linelen
\newdimen\@clnwd
\newdimen\@clnht
\newdimen\@dashdim
\newbox\@dashbox
\newcount\@dashcnt
\catcode`@=12

\catcode`@=11
\font\@linefnt linew10 at 2.4pt
\catcode`@=12

\def\arbreun{\kern-0.4ex
\hbox{\unitlength=.25pt
\picture(60,40)(0,0)
\put(30,0){\droite(0,1){10}}
\put(30,10){\droite(-1,1){30}}
\put(30,10){\droite(1,1){30}}
\put(20,20){\droite(1,1){20}}
\put(10,30){\droite(1,1){10}}
\endpicture}\kern 0.4ex}

\def\arbredeux{\kern-0.4ex
\hbox{\unitlength=.25pt
\picture(60,40)(0,0)
\put(30,0){\droite(0,1){10}}
\put(30,10){\droite(-1,1){30}}
\put(30,10){\droite(1,1){30}}
\put(20,20){\droite(1,1){20}}
\put(30,30){\droite(-1,1){10}}
\endpicture}\kern 0.4ex}

\def\arbretrois{\kern-0.4ex
\hbox{\unitlength=.25pt
\picture(60,40)(0,0)
\put(30,0){\droite(0,1){10}}
\put(30,10){\droite(-1,1){30}}
\put(30,10){\droite(1,1){30}}
\put(50,30){\droite(-1,1){10}}
\put(10,30){\droite(1,1){10}}
\endpicture}\kern 0.4ex}

\def\arbrequatre{\kern-0.4ex
\hbox{\unitlength=.25pt
\picture(60,40)(0,0)
\put(30,0){\droite(0,1){10}}
\put(30,10){\droite(-1,1){30}}
\put(30,10){\droite(1,1){30}}
\put(40,20){\droite(-1,1){20}}
\put(30,30){\droite(1,1){10}}
\endpicture}\kern 0.4ex}

\def\arbrecinq{\kern-0.4ex
\hbox{\unitlength=.25pt
\picture(60,40)(0,0)
\put(30,0){\droite(0,1){10}}
\put(30,10){\droite(-1,1){30}}
\put(30,10){\droite(1,1){30}}
\put(40,20){\droite(-1,1){20}}
\put(50,30){\droite(-1,1){10}}
\endpicture}\kern 0.4ex}

\def\arbreA{\kern-0.4ex
\hbox{\unitlength=.25pt
\picture(60,40)(0,0)
\put(30,0){\droite(0,1){10}}
\put(30,10){\droite(-1,1){30}}
\put(30,10){\droite(1,1){30}}
\endpicture}\kern 0.4ex}

\def\arbreB{\kern-0.4ex
\hbox{\unitlength=.25pt
\picture(60,40)(0,0)
\put(30,0){\droite(0,1){10}}
\put(30,10){\droite(-1,1){30}}
\put(30,10){\droite(1,1){30}}
\put(15,25){\droite(1,1){15}}
\endpicture}\kern 0.4ex}

\def\arbreC{\kern-0.4ex
\hbox{\unitlength=.25pt
\picture(60,40)(0,0)
\put(30,0){\droite(0,1){10}}
\put(30,10){\droite(-1,1){30}}
\put(30,10){\droite(1,1){30}}
\put(45,25){\droite(-1,1){15}}
\endpicture}\kern 0.4ex}

\def\arbreBC{\kern-0.4ex
\hbox{\unitlength=.25pt
\picture(60,40)(0,0)
\put(30,0){\droite(0,1){40}}
\put(30,10){\droite(-1,1){30}}
\put(30,10){\droite(1,1){30}}
\endpicture}\kern 0.4ex}

\def\arbreuut{\kern-0.4ex
\hbox{\unitlength=.25pt
\picture(60,40)(0,0)
\put(30,0){\droite(0,1){10}}
\put(30,10){\droite(-1,1){30}}
\put(30,10){\droite(1,1){30}}
\put(20,20){\droite(0,1){20}}
\put(20,20){\droite(1,1){20}}
\endpicture}\kern 0.4ex}

\def\arbretut{\kern-0.4ex
\hbox{\unitlength=.25pt
\picture(60,40)(0,0)
\put(30,0){\droite(0,1){10}}
\put(30,10){\droite(-1,1){30}}
\put(30,10){\droite(1,1){30}}
\put(30,10){\droite(0,1){20}}
\put(30,30){\droite(1,1){10}}
\put(30,30){\droite(-1,1){10}}
\endpicture}\kern 0.4ex}

\def\arbretuu{\kern-0.4ex
\hbox{\unitlength=.25pt
\picture(60,40)(0,0)
\put(30,0){\droite(0,1){10}}
\put(30,10){\droite(-1,1){30}}
\put(30,10){\droite(1,1){30}}
\put(40,20){\droite(0,1){20}}
\put(40,20){\droite(-1,1){20}}
\endpicture}\kern 0.4ex}

\def\arbreutt{\kern-0.4ex
\hbox{\unitlength=.25pt
\picture(60,40)(0,0)
\put(30,0){\droite(0,1){10}}
\put(30,10){\droite(-1,1){30}}
\put(30,10){\droite(1,1){30}}
\put(10,30){\droite(1,1){10}}
\put(30,10){\droite(0,1){30}}
\endpicture}\kern 0.4ex}

\def\arbrettu{\kern-0.4ex
\hbox{\unitlength=.25pt
\picture(60,40)(0,0)
\put(30,0){\droite(0,1){10}}
\put(30,10){\droite(-1,1){30}}
\put(30,10){\droite(1,1){30}}
\put(50,30){\droite(-1,1){10}}
\put(30,10){\droite(0,1){30}}
\endpicture}\kern 0.4ex}

\def\arbrettt{\kern-0.4ex
\hbox{\unitlength=.25pt
\picture(60,40)(0,0)
\put(30,0){\droite(0,1){10}}
\put(30,10){\droite(-1,1){30}}
\put(30,10){\droite(1,1){30}}
\put(30,10){\droite(-1,2){15}}
\put(30,10){\droite(1,2){15}}
\endpicture}\kern 0.4ex}

\def\arbrebigtut{\kern-0.4ex
\hbox{\unitlength=.25pt
\picture(800,400)(0,0)
\put(300,0){\droite(0,1){100}}
\put(300,100){\droite(-5,3){450}}
\put(300,100){\droite(-1,2){150}}
\put(300,100){\droite(5,3){450}}
\put(300,100){\droite(1,2){150}}
\put(300,100){\droite(0,1){200}}
\put(300,300){\droite(1,1){100}}
\put(300,300){\droite(-1,1){100}}
\put(0,380){$\cdots$}
\put(280,380){$\cdots$}
\put(600,380){$\cdots$}
\put(-160,400){0}
\put(80,420){$r-1$}
\put(200,420){$r$}
\put(340,420){$r+k$}
\put(460,420){$r+k+1$}
\put(740,420){$n$}
\endpicture}\kern 0.4ex}

\def\cube2{\kern-0.4ex
\hbox{\unitlength=.25pt
\picture(800,500)(0,0)
\put(0,0){\droite(1,0){600}}
\put(0,70){\droite(1,0){600}}
\put(0,223){\droite(1,0){600}}
\put(0,0){\droite(2,3){300}}
\put(600,0){\droite(-2,3){300}}
\put(520,-30){\droite(-2,3){300}}
\put(320,-30){\droite(-2,3){300}}
\put(280,140){${\cal C}^2$}
\put(-30,-30){0}
\put(620,-30){1}
\put(290,460){2}
\put(500,-70){${\cal H}_{0,12}$}
\put(300,-70){${\cal H}_{12,0}$}
\put(-100,100){${\cal H}_{2,01}$}
\put(-100,220){${\cal H}_{01,2}$}
\put(340, 230){$123$}
\put(70,230){$303$}
\put(200,30){$321$}
\put(400,30){$141$}
\endpicture}\kern 0.4ex}

\def\associa2{\kern-0.4ex
\hbox{\unitlength=.25pt
\picture(800,500)(0,0)
\put(0,0){\droite(1,0){600}}
\put(0,70){\droite(1,0){600}}
\put(0,223){\droite(1,0){600}}
\put(0,0){\droite(2,3){300}}
\put(80,-30){\droite(2,3){300}}
\put(600,0){\droite(-2,3){300}}
\put(520,-30){\droite(-2,3){300}}
\put(320,-30){\droite(-2,3){300}}
\put(280,140){${\cal K}^2$}
\put(-30,-30){0}
\put(620,-30){1}
\put(290,460){2}
\put(500,-70){${\cal H}_{0,12}$}
\put(300,-70){${\cal H}_{12,0}$}
\put(-100,100){${\cal H}_{2,01}$}
\put(-100,220){${\cal H}_{01,2}$}
\put(380,420){${\cal H}_{1,02}$}
\put(340, 230){$123$}
\put(200,230){$213$}
\put(120,140){$312$}
\put(200,30){$321$}
\put(400,30){$141$}
\endpicture}\kern 0.4ex}

\def\permuto2{\kern-0.4ex
\hbox{\unitlength=.25pt
\picture(800,500)(0,0)
\put(0,0){\droite(1,0){600}}
\put(0,70){\droite(1,0){600}}
\put(0,223){\droite(1,0){600}}
\put(0,0){\droite(2,3){300}}
\put(80,-30){\droite(2,3){300}}
\put(280,-30){\droite(2,3){300}}
\put(600,0){\droite(-2,3){300}}
\put(520,-30){\droite(-2,3){300}}
\put(320,-30){\droite(-2,3){300}}
\put(280,150){${\cal P}^2$}
\put(-30,-30){0}
\put(620,-30){1}
\put(290,460){2}
\put(500,-70){${\cal H}_{0,12}$}
\put(300,-70){${\cal H}_{12,0}$}
\put(-100,100){${\cal H}_{2,01}$}
\put(-100,220){${\cal H}_{01,2}$}
\put(380,420){${\cal H}_{1,02}$}
\put(580,420){${\cal H}_{02,1}$}
\put(340, 230){$123$}
\put(200,230){$213$}
\put(120,140){$312$}
\put(190,40){$321$}
\put(400,140){$132$}
\put(350,40){$231$}
\endpicture}\kern 0.4ex}

\def\arbretroisun{\kern-0.4ex
\hbox{\unitlength=.50pt
\picture(80,100)(0,0)
\put(30,0){\droite(0,1){20}}
\put(30,20){\droite(-1,1){40}}
\put(30,20){\droite(1,1){40}}
\put(60,50){\droite(-1,1){10}}
\put(10,40){\droite(1,1){20}}
\put(70,50){... 1}
\put(70,35){... 2}
\put(70,20){... 3}
\endpicture}\kern 0.4ex}

\def\arbretroisdeux{\kern-0.4ex
\hbox{\unitlength=.50pt
\picture(80,100)(0,0)
\put(30,0){\droite(0,1){20}}
\put(30,20){\droite(-1,1){40}}
\put(30,20){\droite(1,1){40}}
\put(50,40){\droite(-1,1){20}}
\put(0,50){\droite(1,1){10}}
\put(70,50){... 1}
\put(70,35){... 2}
\put(70,20){... 3}
\endpicture}\kern 0.4ex}

\def\doublepeigne{\kern-0.4ex
\hbox{\unitlength=.25pt
\picture(800,240)(0,0)
\put(200,0){\droite(0,01){20}}
\put(200,20){\droite(-1,1){200}}
\put(200,20){\droite(1,1){200}}
\put(20,200){\droite(1,1){20}}
\put(380,200){\droite(-1,1){20}}
\put(90,130){\droite(1,1){100}}
\put(310,130){\droite(-1,1){100}}
\put(80,200){$\cdots$}
\put(280,200){$\cdots$}
\endpicture}\kern 0.4ex}

 \centerline {\grand  Realization of the Stasheff polytope} 

\bigskip \bigskip 
\BB
\centerline {\bf Jean-Louis Loday}
\bigskip \bigskip 
\N {\bf Abstract} {\petit We propose a simple formula for the coordinates of the vertices of the Stasheff polytope (associahedron) and
we compare it to the permutohedron. }

\B

\N {\bf Introduction}. 
\M
The Stasheff polytope $\KKK ^n$, also called associahedron, appeared in the sixties  in the work of Jim Stasheff [St1] on the recognition of loop spaces. It is a
convex polytope of dimension $n$ with one vertex for each planar binary tree with $n+1$ leaves. There are various realizations of $\KKK ^n$ as a polytope in the
literature (cf. [Lee, GKZ, T, St2, D, CFZ]). Here we propose a simple one which has the following advantages, on top of being simple:

--  it  respects the  symmetry, 

-- it  fits with the classical realization of the permutohedron ${\cal P}^n$, 

-- the faces have simple equations.
\M

To any planar binary tree we associate a point in the euclidean space by describing its coordinates
 (which are going to be positive integers) in terms of the
structure of the tree. Explicitly the $i$th coordinate is the product of the number of leaves on the left and on the
right side of the $i$th vertex. The main idea is to start with the permutohedron and to think of it as the truncation of the
standard simplex by some hyperplanes. Truncating only by the {\it admissible} hyperplanes gives the Stasheff polytope.
From the explicit equations of the facets of the permutohedron we compute the coordinates of the intersections of the
admissible hyperplanes and find the result mentioned above.
\B

\N {\bf Convention.} In the euclidean space $\RR^{n}$ the coordinates of a point are denoted $x_1, \cdots , x_n$. We denote by $H$
 the affine hyperplane
 whose equation is: $\sum_{i=1}^{i=n}x_i = {1\over 2}n(n+1)$. We adopt the notation $S(n) =  {1\over 2}n(n+1)$.
\M

\N {\bf 1. A simple realization of the Stasheff polytope}
\M

The {\it Stasheff polytope} $\KKK ^n$ of dimension $n$ (alias associahedron) is a finite cell complex whose $k$-cells are in bijection with the planar trees
having
$n-k+1$ internal vertices and $n+2$ leaves (so it is sometimes denoted $K_{n+2})$, cf. [St1]. 

Let $Y_n$ be the set of planar binary trees with $n+1$ leaves:
$$Y_0 =\{\  |\  \}\ ,\  Y_1= \{\  \arbreA \  \}\ ,\ Y_2=  \{\   \arbreB
,\arbreC \  \}\ ,$$ 
$$ Y_3=  \{\  \arbreun ,\arbredeux ,\arbretrois ,\arbrequatre ,\arbrecinq \ 
\}.$$
The integer $n$ is called the {\it degree} of $t\in Y_n$.
We label the leaves of  $t$ from left to right by $0, 1, \cdots\  $. Then we label the internal vertices  by 
$ 1,2, \cdots\  $. The $i$th vertex is the one which falls in between the leaves $i-1$ and $i$. We denote by $a_i$, resp.
$b_i$, the number of leaves
 on the  left side, resp. right side, of the $i$th vertex.  The product $a_ib_i$ is called the {\it weight} of the $i$th vertex. 
 To the tree $t$ in $Y_{n}$ we associate the point  $M(t)\in
\RR^{n}$ whose $i$th coordinate is the weight of the $i$th vertex:

$$M(t)= ( a_1 b_1, \cdots , a_i b_i, \cdots , a_n b_n) \in \RR^{n}.$$
For instance: $M(\ \arbreA )=(1), \ M(\ \arbreB )= (1,2), \ M(\ \arbreC )=(2,1),$\hfill
$ \ M(\ \arbreun )=(1,2,3), \ M(\ \arbretrois ) =(1,4,1)$.

Observe that the weight of a vertex depends only on the subtree that it determines.

We will show in the next section that all the points   $M(t)\in \RR^{n}$ for $t\in Y_{n}$ lie in the  affine hyperplane $H$. The main point of this paper is
the following result.
\M

\N {\bf 1.1 Theorem.} {\it  The convex hull of the points  $M(t)\in \RR^{n}$, for $t$ a planar binary tree with $n+1$ leaves,  is a realization 
of the Stasheff polytope
 $\KKK ^{n-1}$ (alias associahedron) of dimension $n-1$.
}
\M

The proof will be given in the next section.
\M

Let us recall the definition of the {\it permutohedron} (alias zylchgon). For any permutation $\ss$ in the symmetric
group
$ S_{n}$ acting on the set
 $\{1, \cdots , n\}$ let $M(\ss)\in \RR ^{n}$ be the point with coordinates $M(\ss)=(\ss(1), \cdots , \ss (n))$. By definition the permutohedron
 ${\cal P}^{n-1}$ is
the convex hull  of the ${n}!$ points $M(\ss)$. Observe that the sum of the $n$ coordinates of $M(\ss)$ is $S(n)$, hence all
the
 points
$M(\ss)$ lie in the affine hyperplane $H$.

Under interpreting a permutation as a planar binary tree with levels (cf. for instance [LR]), and then forgetting the levels, one gets a well-defined
 map 
$$\psi : S_{n} \proj Y_{n}.$$
For instance $\psi (1\ 2\ 3)= \arbreun , \psi (1\ 3\ 2)= \psi (2\ 3\ 1)= \arbretrois $.
\M
 
\N {\bf 1.2 Proposition.} {\it  The Stasheff polytope  $\KKK ^{n-1}$ as defined above contains the permutohedron 
 ${\cal P}^{n-1}$ and, for $\ss\in S_{n}$ and $t\in Y_{n}$,  the following are equivalent:

(a) $M(\ss) = M(t)$,

(b) $\psi\mm (t) = \{\ss\}$,

(c) $\psi(\ss) = t$ and for all $i$ either $a_i=1$ or $b_i=1$,

(d)  $\psi(\ss) = t$ and  the permutation $\ss$ has the following property: either $\ss(1)=n$ or $\ss (n)=n$ and the remaining permutation $\ss'\in
S_{n-1}$, obtained by deleting $n$,  satisfies the same property (ad libitum). }
\M

The proof will be given in the next section. Observe that, as a consequence, there are $2^{n-1}$ points which are both
vertices of  $\KKK ^{n-1}$ and of ${\cal P}^{n-1}$.
\M

\N {\bf Example.} For $n=3$, the common vertices are the points corresponding to the permutations $(1\ 2\ 3), (2\ 1\ 3),(3\ 1\ 2), (3\ 2\ 1)$.
  For $n=4$ we get  $(1\ 2\ 3\ 4), (2\ 1\ 3\ 4),(3\ 1\ 2\ 4), (3\ 2\ 1\ 4), 
(4\ 1\ 2\ 3), (4\ 2\ 1\ 3),(4\ 3\ 1\ 2), (4\ 3\ 2\ 1)$.
\B

\N {\bf 2. Equation for the facets of the permutohedron and of the Stasheff polytope}
\M
The permutohedron can be obtained
by truncating the standard simplex along some hyperplanes, one per each cell of $\DD^{n-1}$ (except the big cell).
Truncating only along ``admissible" ones gives the Stasheff polytope (cf. [St2]). We give the explicit equations of these hyperplanes.
As a consequence we get the results announced in the first section.
\M

 \N {\bf 2.1 Shuffles and hyperplanes.}  The intersection of the hyperplane $H$ with the quadrant $\{(x_1, \cdots , x_n) \mid x_i\geq 0 \hbox{ for all }
i\}$ is (homothetic to) the standard
 simplex  $\DD^{n-1}$. We still call this intersection the standard simplex.

 The $(k-1)$-cells of $\DD^{n-1}$, for $k=1, \cdots n-1$, are indexed by the $(k,n-k)$-shuffles of $(1,\cdots , n)$, i.e. the partitions $\omega =
(\oo_1 \cdots \oo_k \vert \oo_{k+1} \cdots
\oo_n)$
 of $\{1,\cdots , n\}$ into 2 nonempty subsets. Two partitions $\oo$ and $\oo'$ are the same if they differ only by the order of the
integers in each subset. It is sometimes necessary to take a representative, in which case we will assume that
$\oo_1 < \cdots <\oo_k$ and $ \oo_{k+1}< \cdots < \oo_n$. We denote by
$\bar {\oo}$ the dual cell, i.e $\bar {\oo}= ( \oo_{k+1} \cdots \oo_n \vert \oo_1 \cdots \oo_k)$.  

We  associate to the shuffle $\omega = (\oo_1 \cdots \oo_k\vert  \oo_{k+1} \cdots \oo_n)$ the polynomial
$$\displaylines{
p_{\oo}\row x1n :=  \hfill \cr
\hfill (n-k)(x_{\oo_1}+ \cdots + x_{\oo_k}) - k(x_{\oo_{k+1}}+ \cdots + x_{\oo_n}) +{1\over 2}nk(n-k),\cr
}$$
and the hyperplane  ${\cal H}_{\oo}$ defined by $p_{\oo}\row x1n =  0$. Observe that ${\cal H}_{\oo}$ and ${\cal H}_{\bar
{\oo}}$ are parallel.
\M

\N {\bf 2.2 Recall on the permutohedron.} For any $\ss\in S_n$ the point $M(\ss)$ lies in the affine hyperplane $H$. The
$k$-cells of the permutohedron
${\cal P}^{n-1}$ can be indexed by the partitions of $\{1, \cdots , n\}$ into $n-k$ subsets. In particular the facet 
indexed by the shuffle $\oo$ lies in the hyperplane ${\cal H}_{\oo}$. The $0$-cells, that is the vertices $M(\ss)$, are
indexed by the partitions $\oo(\ss)$ with $n$ subsets. The relationship with the permutations is given by
$$\oo(\ss) = \ss\mm.$$
 The vertex $M(\ss)=\row x1n$ lies in the hyperplane ${\cal H}_{\oo}$ if and
only if the partition  $\oo(\ss)$ is a refinement of $\oo$. Indeed, if it is so, then $\{ x_{\oo_1}, \cdots , x_{\oo_k}\}= \{ 1, \cdots,
k\}$ and $ x_{\oo_1}+
\cdots + x_{\oo_k}= S(k)$. It also implies $x_{\oo_{k+1}}+ \cdots +x_{\oo_n}= S(n) - S(k)$ and therefore $ p_{\oo}(M) =0$. Observe that if $\oo(\ss)$ is not a
refinement of $\oo$, then $ p_{\oo}(M) >0$.

In conclusion, the permutohedron ${\cal P}^{n-1}$ is the truncation of the standard simplex lying in $H$ by the hyperplanes
${\cal H}_{\oo}$ for all shuffles $\oo$.
\M
\N {\bf Example: $n=3$.}

$$\permuto2$$
\BB
$$\displaylines{
M(1\vert 2\vert 3) = (1,2,3), \ M(2\vert 1 \vert 3) = (2,1,3), \ M(2\vert 3\vert 1) = (3,1,2), \ \cr
M(3\vert 2 \vert 1) = (3,2,1), \ M(3\vert 1 \vert 2) = (2,3,1), \ M(1\vert 3\vert 2) = (1,3,2). \ \cr
}$$
\M

\N {\bf 2.3 Planar trees and admissible shuffles.} We denote by $T_n$ the set of planar trees with $n+1$ leaves, $n\geq
0$ (and one root) such that the valence of each internal vertex is at least 2.  Here are the first of them:
$$T_0 = \{\vert\} ,\qquad \ T_1 = \{\ \arbreA\} ,\qquad T_2 = \{\ \arbreB , \arbreC ; \arbreBC\}$$ 
$$T_3 = \{\ \arbreun, \arbredeux, \arbretrois, \arbrequatre, \arbrecinq; \arbreuut, \arbretut, \arbretuu, \arbreutt, \arbrettu; \arbrettt \}. $$
The integer $n$ is called the {\it degree} of $t\in T_n$. The set $T_n$ is the disjoint union of the
sets $T_{n,k}$  made of the planar trees which have $n-k+1$ internal vertices. For instance $T_{n,1}= Y_n $ since it is made of the
planar binary trees. On the other extreme the set $T_{n,n}$ has only one
 element, which is the planar tree with one vertex. It is sometimes called a {\it corolla}.  So we have 
$$T_n = T_{n,1} \cup \cdots \cup T_{n,n}.$$

A tree $t$ is called a refinement of the tree $t'$ if $t'$ can be obtained from $t$ by contracting to a point some of the internal edges. So any tree is a refinement of the
corolla. Any tree in $T_{n, n-1}$ is of the form 
 
\B
$$y(r,k) = \arbrebigtut$$
\B

Observe that, once $n$ is fixed,  it is completely determined by $r$ and $k$, whence the notation $y(r,k) $. We associate to it the 
shuffle $\oo(y(r,k))= ( r+1\ r+2\ \cdots \ r+k \vert \ \cdots \ )$.

 A shuffle $\omega = (\oo_1 \cdots \oo_k \vert \oo_{k+1} \cdots \oo_n)$ such that
$\oo_1 < \cdots <\oo_k$ and $ \oo_{k+1}< \cdots < \oo_n$ is called {\it admissible} if the first part is a sequence of {\it consecutive} integers, i.e. if $\oo$ is of the
form $(r+1\ r+2 \cdots \ r+k \vert \ \cdots \ )$. Observe that there is a bijection between  the admissible shuffles and the planar tree with  2 internal
vertices. We denote by $\oo(t)$ the partition associated to $t$.
\M

\N {\bf 2.4 Recall on the Stasheff polytope.} It is shown in [St2, Appendix] that the Stasheff polytope can be obtained from the standard simplex by truncating
along the hyperplanes corresponding to the admissible shuffles. We will show that the points $M(t)$ defined in section 1 are indeed the vertices of the
polytope in $H$ defined by the equations $ p_{\oo}(M)\geq 0 $ for $\oo$ an admissible shuffle.
\M

\N {\bf  2.5 Lemma.} {\it  For any tree $t\in Y_{n}$ the coordinates of the point $M(t)=(x_1, \cdots , x_n) \in \RR^{n}$ satisfy the relation
 $\sum _{i=1}^{i=n} x_i = {1\over 2}n(n+1)$. Hence one has $M(t)\in H$.}
\M

\N {\it Proof.}  Any planar binary tree  $t$ in $Y_{n}, n\geq 1,$ is the grafting of its left part $t^l$ and its right part $t^r$
(cf. [LR]), so
$t = t^l\vee t^r$. For instance $\arbreA = \vert \vee \vert\ $ . It follows 
 from the definition of $M(t)$ that
$$M(t)= (M(t^l),pq, M(t^r))$$
where $p$ is the number of leaves of $t^l$ and  $q$ is the number of leaves of $t^r$ (we initialize with $M(\vert) =
\emptyset$). For
$M = (x_1,\cdots, x_n)\in
\RR^{n}$, let
 $\phi(M):=x_1+\cdots + x_n$. We have  $\phi(M(\arbreA))= 1$ since $M(\arbreA)= (1)$.
Suppose, by induction, that, $\phi(M(t))= {1\over 2}k(k+1)$ for any $t\in Y_k$, $k<n$. Then, for $t\in Y_{n}$, we get
 $\phi(M(t))= \phi(M(t^l)) + pq + \phi(M(t^r))$. By induction $\phi(M(t^l))=  {1\over 2}(p-1)p$ and 
$\phi(M(t^r))=  {1\over 2}(q-1)q$. Hence we get
$$\eqalign{
\phi(M(t))&= {1\over 2}(p-1)p+ pq + {1\over 2}(q-1)q\cr
          &= {1\over 2}(p^2 +2 pq + q^2 -p -q)\cr
          &= {1\over 2}(n+1)n\cr
}$$
since $n+1 = p+q$. Hence we have proved that $M(t)$ belongs to the affine hyperplane $H$ of $\RR^{n}$.
\hfill $\square$
\M

\N {\bf  2.6 Proposition.} {\it  Let    $\oo$ be an admissible shuffle. For any tree $t\in Y_{n}$ the point $M(t)$ lies in ${\cal H}_{\oo}$ if and
only if the partition $\oo(t)$ is a refinement of $\oo$. If not, then $ p_{\oo}(M(t)) >0$.}
\M

\N {\it Proof.} Let $y(r,k) $ be the tree corresponding to $\oo$. When  $\oo(t) $ is a refinement of $\oo$,  $t$ is a planar
 binary tree  such that
 by contracting some internal edges we can obtain the tree $y(r,k)$ with one internal
edge. Let 
$(x_1, \cdots , x_n)$ be the coordinates of $M(t)$. From the structure of $t$ it follows that the subtree which contains
 the leaves number $r$ to $r+k$
 is a
tree of degree $k$ (its root becomes the only internal edge of $y(r,k)$). Hence we get $ x_{r+1}  + x_{r+2}  + \cdots +
x_{r+k}= S(k)$. Since 
$\sum_{i=1}^{i=n}x_i =S(n)$, we get $\sum_{i=1}^{i=r}x_i + \sum_{i=r+k+1}^{i=n}x_i =S(n)-S(k)$. 
\S
We already know that $M(t)\in H$ by Lemma 2.5. The conclusion $M(t)\in H_{\oo}$ follows from the following computation:
$$\eqalign{
p_{\oo}(M) &= (n-k) S(k) -k(S(n) - S(k)) + {1\over 2}nk(n-k)\cr
&=n S(k) - k S(n) + {1\over 2}nk(n-k)\cr
&= {1\over 2}nk(k+1-n-1+n-k)\cr
&=0.\cr
}$$
If $\oo(t)$ is not a refinement of $\oo$, then one of the values of $ x_{r+1}  + x_{r+2}  + \cdots +
x_{r+k}$ is at least $1\times (k+1)= k+1$ and the sum
 of the others
 is at least $S(k-1)$. Hence $\sum_{i=r+1}^{i=r+k}x_i \geq S(k-1) + k+1 > S(k)$. Since $\sum_{i=1}^{i=n}x_i =S(n)$, we get
$ p_{\oo}(M(t)) >0$.
\hfill $\square$
\M

\N {\it Proof of Theorem 1.1.} Since the Stasheff polytope is the truncation of the standard simplex by the hyperplanes corresponding to the admissible
shuffles, the vertex corresponding to the tree $t$ is the intersection of the hyperplanes $H, {\cal H}_{\oo_1}, \cdots ,  {\cal
H}_{\oo_{n-1}}$, where 
${\oo_1}, \cdots , {\oo_{n-1}}$ are the admissible shuffles corresponding to the trees with 2 vertices which admit $t$
as a refinement. In order to show that the convex hull of the points $M(t)$ is a realization of the Stasheff polytope, it is
sufficient to show that 
$$H \cap {\cal H}_{\oo_1}\cap \cdots \cap  {\cal H}_{\oo_{n-1}}= \{ M(t) \}.$$
Since we know that this intersection is a point, it is sufficient to prove that $M(t)$ lies in each hyperplane. For $H$ this is Lemma 2.5. 
For ${\cal H}_{\oo_i}$ this is Proposition 2.6 and we are done.
\hfill $\square$
\M

\N {\bf Example: $n=3$.}

$$\associa2$$
\BB
$$\displaylines{
M(\ \arbreun) = (1,2,3), \ M(\ \arbredeux) = (2,1,3), \ M(\ \arbrequatre) = (3,1,2), \ \cr
M(\ \arbrecinq) = (3,2,1), \ M(\ \arbretrois) = (1,4,1). \ \cr
}$$
\M
\N {\it Proof of Proposition 1.2.} Recall that the surjective map $\psi : S_{n} \proj Y_{n}$ is defined by using the interpretation of the permutations as
planar binary trees {\it with levels}, cf. [L-R].  In a tree with levels each internal vertex has a level (ranging from 1 to $n$), and there is only one vertex per
level. The permutation is obtained by taking $\ss(i) =$ the level of the $i$th vertex. For instance
$$\arbretroisun\qquad \qquad \arbretroisdeux$$

\N (a) $\Rightarrow$ (b)  Since all the points $M(\ss)$ are distinct, the equality $M(\ss)=M(t)$ implies that, for a fixed $t$,  there is only one $\ss$ such that
$\psi(\ss) = t$.

\N (b) $\Rightarrow$ (c) The equality $\psi(\ss) = t$ is clear from the assumption. If some internal vertex of $t$ has several leaves on both sides, then there
are several ways of lifting it as a leveled tree.

\N (c) $\Rightarrow$ (d) The root vertex has either $a_i=1$ or $b_i=1$. In the first case, it means that $i=1$ and the
permutation is of the form $( n\ \cdots
\ )$. In the second case it means that $i=n$ and the permutation is of the form $(\ \cdots\ n)$. And so on.

\N (d) $\Rightarrow$ (a) Let $\ss$ be a permutation as described in condition (d). Then the $i$th vertex of the associated tree has one leaf on one side and has
$\ss(i)$ ($=$ the degree of the subtree it generates) leaves on the other side. Hence the weight of the $i$th vertex is $1\times \ss(i) = \ss(i)$. It follows that
$M(\ss)=M(t)$. \hfill $\square$
\M

\N {\bf 2.7 The cube.} There is a family of cubes ${\cal C}^n$ which fits nicely with the families of Stasheff polytopes
and permutohedrons. It is defined as follows.

The convex polytope ${\cal C}^{n-1}$ is defined in $H$ by the equations $ p_{\oo}(M)\geq 0 $ for  $\oo = (1\ 2 \ \cdots \ i
\vert i+1\ \cdots \ n)$ and 
$\oo = (i+1\ \cdots \ n \vert 1\ 2 \ \cdots \ i)$, $i$ ranging from 0 to $n-1$. So the facets of ${\cal C}^{n-1}$ are in the
hyperplanes
$$ {\cal H}_{1, 2 \cdots \ n}\ ,\  {\cal H}_{12, 3 \cdots \ n}\ ,\ \cdots \ ,  {\cal H}_{1 \cdots \ n-1, n}\ ,$$
and 
$$ {\cal H}_{2 \cdots \ n,1 }\ ,\  {\cal H}_{3 \cdots \ n,12 }\ ,\ \cdots \ ,  {\cal H}_{n, 1 \cdots \ n-1}\ .$$

Let us introduce the set $Q_{n} := \{0,1\}^{n-1}$ and the map $\phi : Y_{n} \proj Q_{n}$ given by 
$\phi(t) = (\ee_1, \cdots , \ee_{n-1})$ where $\ee_i = 0$ (resp. $\ee_i = 1$) if the $i$th leaf of $t$ is pointing to the left
(resp. to the right) (cf. [LR]). We code the vertices of the cube by the elements of $Q_{n}$. With the definition given above
for ${\cal C}^{n-1}$, the coordinates of the vertex $M(\ee)=\row x1n$ for 
$\ee = \row {\ee}1{n-1} \in Q_{n}$ are 
$$x_i = i - \ee_{i-1}(n-i+1)(i-1)+\ee_{i}(n-i)(i), \quad \hbox {for } i=1, \cdots , n.$$
\M

\N {\bf Example: $n=3$.}

$$\cube2$$
\BB
One can check that ${\cal C}^{n-1}$ contains the Stasheff polytope ${\cal K}^{n-1}$ and hence the permutohedron ${\cal
P}^{n-1}$. The only common vertices of ${\cal C}^{n-1}$ and ${\cal P}^{n-1}$ are the two points $(1, 2, \cdots , n)$ and  $(n ,
n-1,
\cdots , 1)$. The polytopes ${\cal C}^n$ and ${\cal K}^n$ have $n$ vertices in common. They are characterized by the
following equivalent conditions, where $t\in Y_{n}$ and $\ee\in Q_{n}$ :
\M

(a) $M(t) = M(\ee)$,

(b) $\phi\mm (\ee) = \{t\}$,

(c) $\phi(t) = \ee$ and $\ee$  is of the form $(1, \cdots , 1,0, \cdots , 0)$,

(d)  $\phi(t) = \ee$ and   $t$ is of the form:
$$\doublepeigne$$
We remark that these common vertices are on the shortest path (for the Tamari order, see below) in ${\cal K}^n$ from the minimal vertex $(1, 2, \cdots , n)$ to
the maximal vertex $(n , n-1, \cdots , 1)$.
\M

\N {\bf 2.8 Relationship with the poset structure.} The three sets $S_{n}$, $Y_{n}$ and $Q_{n}$ can be equipped with a poset structure such that the
maps $\psi$ and $\phi$ are maps of posets, cf. [LR]. On $S_{n}$ it is called the {\it weak Bruhat order}, on $Y_{n}$ it is called the {\it Tamari order}, and
on $Q_{n}$ it is the lexicographic order. These orders are induced by the following {\it covering} relations.

Let $s_i$ be the permutation which exchanges $i$ and $i+1$. Then $\{s_1, \cdots, s_{n-1}\}$ is a set of generators of
$S_{n}$. For two permutations $\ss$ and
$\ss'$, $\ss < \ss'$ is said to be a covering relation for the weak Bruhat order if $\ss' = s_i \ss $ for some $i$ and the length of $\ss'$ is greater than the
length of $\ss$ (recall that the length of $\ss$ is the minimal number of generators necessary to write $\ss$ in terms of the $s_i$'s).
 So, for $S_3$ the poset structure is: 
$$\matrix{
 & &1& & \cr
& & & &\cr
 & \swarrow & &\searrow & \cr
 s_1& & & & s_2 \cr
& & & &\cr
 \downarrow& & & &\downarrow \cr
& & & &\cr
 s_2s_1 & & & &s_1 s_2  \cr
 &\searrow & &\swarrow & \cr
& & & &\cr
 & & s_1s_2  s_1=s_2  s_1s_2  & & \cr
}$$
 For two planar binary trees $t$ and
$t'$, $t < t'$ is said to be a covering relation for the Tamari order if $t'$ can be obtained from $t$ by changing locally the pattern $\arbreB$ to the pattern
$\arbreC$. So for $Y_3$ the poset structure is: 
$$\matrix{
 & &\arbreun& & \cr
 & \swarrow & & & \cr
\arbredeux& & &\searrow & \cr
& & & &\cr
 \downarrow& & & &\arbretrois \cr
& & & &\cr
\arbrequatre & & & \swarrow& \cr
 &\searrow & & & \cr
 & &\arbrecinq & & \cr
}$$

For $Q_{n}$ the lexicographic order is induced by the following covering relation: $\ee < \ee'$ if and only if the values of $\ee$ and $\ee'$ are the same
except at one place, let say the $i$th one, for which $\ee_i = 0$ and $\ee'_i = 1$.
 So for $Q_3$ the poset structure is: 
$$\matrix{
 & &(0,0)&  \cr
 & \swarrow & &\searrow   \cr
& & & &\cr
(1,0)& & & &(0,1) \cr
& & & &\cr
 &\searrow & &\swarrow   \cr
 & &(1,1) &  \cr
}$$

\N {\bf 2.9 Proposition.} {\it The realization of the polytopes ${\cal P}^{n-1}$, ${\cal K}^{n-1}$ and ${\cal C}^{n-1}$
described above are such that all the edges are precisely the covering relations of the posets.

The projection $\pi$ of the vertices on the oriented axis $NS$ where $N=(1, 2, \cdots , n)$ and $S=(n , n-1, \cdots , 1)$
respects the order.}
\hfill $\square$
\M

\N {\bf 2.10 Length of the edges.} It is immediate to verify that the length of the edges in the permutohedron case is
$\sqrt {2}$ and in the case of the cube $(n-i)i\sqrt{2}$ for $i=1, \cdots n-1$. In the case of the Stasheff polytope it is also an
integral multiple of
$\sqrt {2}$. Indeed if the covering relation from $t$ to $t'$ moves the $i$th vertex of $t$ to the $j$th vertex of $t'$, then
the corresponding edge is of length $a_ib_j \sqrt{2}$. So the maximum length of an edge in ${\cal K}^{n}$ is
$m(m+1)\sqrt{2}$ if
$n=2m$ and $m^2\sqrt{2}$ if $n=2m-1$. 
\M
\N {\bf 2.11 Barycenter.} F. Chapoton observed that the barycenter $B$ of the vertices of ${\cal K}^{n-1}$ is 
$({n+1\over 2} , \ldots , {n+1\over 2} )$. So the permutohedron, the Stasheff polytope and the ``cube" have the same barycenter. 
As a consequence, for any vertex $M$ of the Stasheff polytope one has
$$ \overrightarrow{BM} = \sum_i \overrightarrow{BP_i}$$
where the points $P_i$ are the vertices of the permutohedron corresponding to the permutations whose image by $\psi$
 is the tree corresponding to $M$ (cf. Proposition 1.2).

Similarly, if $N$ is a vertex of the ``cube", then 
$$ \overrightarrow{BN} = \sum_j \overrightarrow{BM_j}$$
where the points $M_j$ are the vertices of the Stasheff polytope related to $N$ by $\phi$ (cf. 2.7).
\vfill
\eject

A projection of the Stasheff polytope ${\cal K}^3$ with the coordinates of its vertices in ${\bf R}^4$.
\BB
\epsfbox {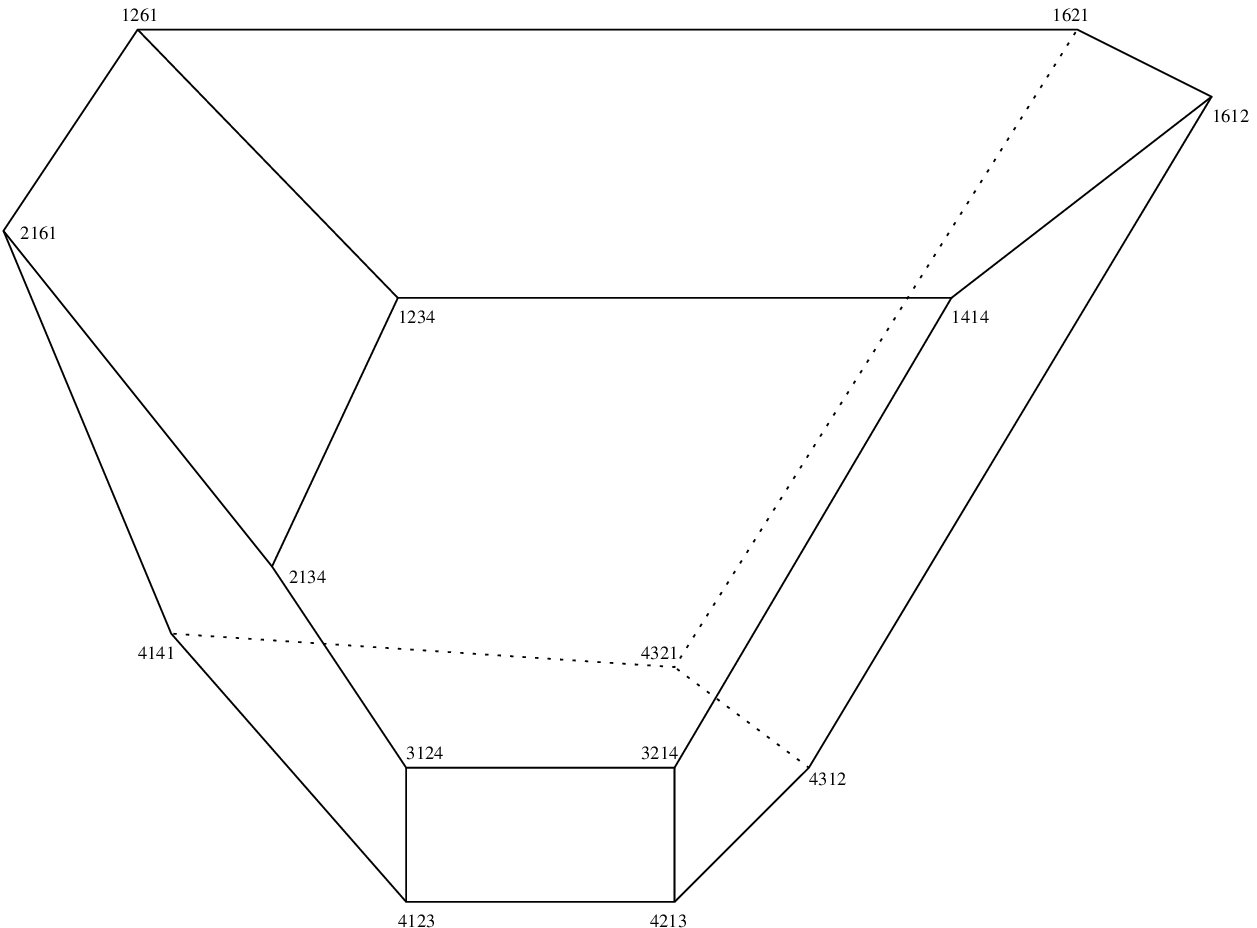}

\BB

\vfill 
\eject

\centerline {\bf References}
\B

\N [CFZ]  Chapoton, F.;  Fomin, S.; Zelevinsky, A. {\it Polytopal realizations of generalized associahedra}. Bulletin Canadien de Math\'ematiques 
(to appear) {\tt [ArXiv: math.CO/0202004]}
\S 
\N [D] Devadoss, S. L. {\it Tessellations of moduli spaces and the mosaic operad.} Homotopy
invariant algebraic structures (Baltimore, MD, 1998), 91--114, Contemp. Math., 239, Amer. Math. Soc.,
Providence, RI, 1999.
\S
\N [GKZ] Gelfand, I. M.; Kapranov, M. M.; Zelevinsky, A. V. Discriminants, resultants, and multidimensional determinants. Mathematics: Theory and
Applications. Birkh\" auser Boston, Inc., Boston, MA, 1994. x+523 pp. 
\S
\N [Lee] Lee, C. W. {\it 
The associahedron and triangulations of the $n$-gon}. 
European J. Combin. 10 (1989), no. 6, 551--560. 
\S
\N [LR] Loday, J.-L., and Ronco M.,  {\it Order structure on the algebra of permutations and of planar binary trees,} J. Algebraic Combinatorics, 15(3) (2002),
253--270.
\S
\N [St1] Stasheff, J. D. {\it Homotopy associativity of $H$-spaces}. I, II. Trans. Amer. Math. Soc. 108 (1963), 275-292;
 ibid. 293--312. 
\S
\N [St2] Stasheff, J. D.  {\it From operads to ``physically" inspired theories.} Operads: Proceedings of
Renaissance Conferences (Hartford, CT/Luminy, 1995), 53--81, Contemp. Math., 202, Amer. Math. Soc.,
Providence, RI, 1997.
\S
\N [T] Tonks, A. {\it Relating the associahedron and the permutohedron.} Operads: Proceedings of
Renaissance Conferences (Hartford, CT/Luminy, 1995), 33--36, Contemp. Math., 202, Amer. Math. Soc.,
Providence, RI, 1997. 
\BB

  Institut de Recherche Math\'ematique Avanc\'ee,

    CNRS et Universit\'e Louis Pasteur

    7 rue R. Descartes,

    67084 Strasbourg Cedex, France

    Courriel : loday@math.u-strasbg.fr
\M

\hfill 9 December 2002
\end

\N {\bf  Lemma.} {\it  For any tree $t\in Y_{n}$ the coordinates of the point $M(t)=(x_1, \cdots , x_n) \in \RR^{n}$ satisfy the relation
 $\sum _{i=0}^{i=n} x_i = {1\over 2}n(n+1)$. Hence one has $M(t)\in H$.}
\M

\N {\it Proof.}  Any planar binary tree  $t$ in $Y_{n}$ is the grafting of its left part $t^l$ and its right part $t^r$, so $t = t^l\vee t^r$. It is clear
 from the construction of $M(t)$ that
$$M(t)= (M(t^l),pq, M(t^r))$$
where $p$ is the number of leaves of $t^l$ and  $q$ is the number of leaves of $t^r$. For $M = (x_1,\cdots, x_n)\in \RR^{n}$, let
 $\phi(M):=x_1+\cdots + x_n$. We have  $\phi(M(\arbreA))= 1$ since $M(\arbreA)= (1)$.
Suppose, by induction, that, $\phi(M(t))= {1\over 2}k(k+1)$ for any $t\in Y_k$, $k\leq n$. Then, for $t\in Y_{n}$, we get
 $\phi(M(t))= \phi(M(t^l)) + pq + \phi(M(t^r))$. By induction $\phi(M(t^l))=  {1\over 2}(p-1)p$ and 
$\phi(M(t^r))=  {1\over 2}(q-1)q$. Hence we get
$$\eqalign{
\phi(M(t))&= {1\over 2}(p-1)p+ pq + {1\over 2}(q-1)q\cr
          &= {1\over 2}(p^2 +2 pq + q^2 -p -q)\cr
          &= {1\over 2}(n+1)n\cr
}$$
since $n+1 = p+q$. Hence we have proved that $M(t)$ belongs to the affine hyperplane $H$ of $\RR^{n}$.
\hfill $\square$
\M